\documentclass[11pt]{article}
\usepackage{amscd}
\usepackage{amsfonts}
\usepackage{amsmath}
\usepackage{booktabs}
\usepackage{url}
\usepackage{hyperref}

\newtheorem{theorem}{Theorem}
\newtheorem{example}{Example}
\newtheorem{lemma}{Lemma}

\newcommand{\F}{{\mathbb F}}
\newcommand{\B}{\mathcal{B}}

\newcommand{\gauss}[3]{\genfrac{[}{]}{0pt}{}{#1}{#2}_{#3}}
\newcommand{\subspaces}[3]{\linearlattice_{#1,#2}(#3)}
\newcommand{\gl}[2]{\generallinear_{#1}(#2)}
\newcommand{\ns}[2]{\normalizersinger_{#1}(#2)}

\DeclareMathOperator\wt{wt}

\DeclareMathOperator\linearlattice{L}
\DeclareMathOperator\generallinear{GL}
\DeclareMathOperator\normalizersinger{N}

\begin{document}

% title
\title{$q$-Analogs of Packing Designs}
\author{
{\sc Michael Braun}\footnote{The authors are with the University of Applied Sciences Darmstadt, Germany, corresponding author's email: \url{michael.braun@h-da.de}}
\and 
{\sc Jan Reichelt}}

\maketitle

% abstract
\begin{abstract}
A $P_q(t,k,n)$ $q$-packing design is a selection of $k$-subspaces of $\F_q^n$ such that each $t$-subspace is contained in at most one element of the collection. A successful approach adopted from the Kramer-Mesner-method of prescribing a group of automorphisms was applied by Kohnert and Kurz to construct some constant dimension codes with moderate parameters which arise by $q$-packing designs. In this paper we recall this approach and give a version of the Kramer-Mesner-method breaking the condition that the whole $q$-packing design must admit the prescribed group of automorphisms. Afterwards, we describe the basic idea of an algorithm to tackle the integer linear optimization problems representing the $q$-packing design construction by means of a metaheuristic approach. Finally, we give some improvements on the size of $P_2(2,3,n)$ $q$-packing designs.
\end{abstract}

\newpage

\section{Introduction}\label{sec:intro}
A \emph{$q$-packing design} $\B$ with parameters $P_q(t,k,n)$, also called \emph{$q$-analog of a packing design}, is a collection of $k$-subspaces of an $n$-dimensional vector space $\F_q^n$ over the finite field with $q$ elements,
\[
\B\subseteq \subspaces{n}{k}{q}:=\{K\le \F_q^n\mid \dim(K)=k\},
\]
such that each $t$-subspace of $\F_q^n$ is contained in at most one element of $\B$:
\[
\forall\, T\in \subspaces{n}{t}{q}: |\{K\in \B\mid T\subseteq K\}|\le 1.
\]

An equivalent characterization is the following one:

\begin{lemma}\label{lem:alternative}
A set $\B\subseteq\subspaces{n}{k}{q}$ is a $P_q(t,k,n)$ $q$-packing design if and only if the intersection of every two elements of $\B$ is at most $(t-1)$-dimensional:
\[
\forall\, K,K'\in\B: \dim(K\cap K')\le t-1.
\]
\end{lemma}

The maximal possible cardinality of a $P_q(t,k,n)$ $q$-packing design will be denoted by $B_q(t,k,n)$---it is called the \emph{$q$-packing number}. If the expression
\[
\gauss{n}{k}{q}:=\genfrac{}{}{}{}{(q^n-1)(q^{n-1}-1)\cdots(q^{n-k+1}-1)}{(q^k-1)(q^{k-1}-1)\cdots(q-1)}
\]
denotes the \emph{$q$-Binomial coefficient} (also called \emph{Gaussian number}), the $q$-packing number is bounded by the \emph{packing bound}
\[
B_q(t,k,n)\le \genfrac{}{}{}{}{\gauss{n}{t}{q}}{\gauss{k}{t}{q}}.
\]
This bound is reached if and only if each $t$-subspace is contained in \emph{exactly} one element of the collection. In this case the $q$-packing design is called a \emph{$q$-Steiner system} with parameters $S_q(t,k,n)$. The long standing issue if non-trivial cases for $S_q(t,k,n)$, $t>1$, do exist at all was recently solved by constructing the first $S_2(2,3,13)$ $q$-Steiner systems \cite{BEO+12}.

The interest in the construction of $q$-packing designs of maximal size have increased since Koetter and Kschischang published their trend-setting paper~\cite{KK08} to an algebraic approach to random network codes. A \emph{random network code} $\mathcal{C}$ is a collection of subspaces of $\F_q^n$. The corresponding metric between subspaces of $\F_q^n$ is given by the \emph{subspace distance} 
\begin{align*}
d(S,T)&=\dim(S+T)-\dim(S\cap T)\\
&= \dim(S)+\dim(T) - 2\dim(S\cap T).
\end{align*} 

If all subspaces of $\mathcal{C}$ have the same dimension $k$ we speak of \emph{constant dimension codes} and call $\mathcal{C}$ an $[n,k,d,s]_q$ code if and only if in addition $\mathcal{C}$ has the minimum distance $d$. Constant dimension codes with parameters $[n,k,2(k-t+1),s]_q$ are related to $P_q(t,k,n)$ $q$-packing designs with cardinality $s$---in fact they are the same~\cite{KK08a}.

A table of lower bounds on $B_q(t,k,n)$ is given in~\cite{KK08a} where the authors list a table for the parameters $q=2$, $t=2$, $k=3$, and  $6\le n\le 14$ (see Table~\ref{tab:old}). Most values are still valid. The upper bound is given by the packing bound.

\begin{table}[!htbp]
\centering
\caption{Lower and upper bounds on $B_2(2,3,n)$}\label{tab:old}
\begin{tabular}{lrrr}
\toprule
$n$ & Lower & Upper & Reference\\
\midrule
$6$ & $77$ & $93$ & \cite{KK08a}\\
$7$ & $304$ & $381$ & \cite{KK08a}\\
$8$ & $1275$ & $1542$ & \cite{KK08a}\\
$9$ & $5694$ & $6205$ & \cite{EV08b}\\
$10$& $21483$ & $24893$ & \cite{KK08a}\\
$11$& $79833$ & $99718$ & \cite{KK08a}\\
$12$& $315315$ & $399165$ & \cite{KK08a}\\
$13$&Ê$1597245$ & $1597245$ & \cite{BEO+12}\\
$14$&Ê$4177665$ & $6390150$ & \cite{KK08a}\\
\bottomrule
\end{tabular}
\end{table}

Kohnert and Kurz~\cite{KK08a} used a well-known approach for the generation of their $q$-packing designs---the so-called \emph{Kramer-Mesner-method}~\cite{KM76}. The basic idea of that approach is the representation of the construction problem as a Diophantine linear system of equations (in case of e.\,g. designs, designs over finite fields), respectively as an integer linear optimization problem (in case of e.\,g. linear codes, arcs, blocking sets, constant dimension codes, $q$-packing designs). In general, the parameters of interest produce large optimization problems. In order to obtain a reduction of the optimization task the search will be restricted to objects that admit a certain group of automorphisms.

In this paper we recall the Kramer-Mesner-method briefly for the case of $q$-packing designs and describe a refined approach of the Kramer-Mesner-method for integer linear programming by breaking the condition that the whole object must admit the prescribed group of automorphisms. Then, we introduce an approach to tackle the integer linear optimization problem defined by the $q$-packing design construction. Finally, we give some improvements of the lower bound on $B_2(2,3,n)$ for some values $n$.

\section{Recalling the Kramer-Mesner-Approach}\label{sec:km}

The \emph{incidence matrix} $A_{t,k}=(a_{TK})$ between all $t$- and all $k$-subspaces of $\F_q^n$ is defined by
\[
a_{TK}:=\left\{\begin{array}{cl}
1, &\text{if $T\subseteq K$}\\
0, &\text{otherwise.}
\end{array}\right.
\] 

The construction of $q$-packing designs can be formulated as a selection of columns of the incidence matrix:

\begin{lemma}[\cite{KK08a}]
The set of $P_q(t,k,n)$ $q$-packing designs with cardinality $s$ bijectively corresponds to the set of 0-1-vectors $x$ such that $A_{t,k}x$ is a 0-1-vector with weight $\wt(x) = s$. If $x$ denotes such a 0-1-vector the corresponding $q$-packing design is given by $\B=\{K\mid x_K=1\}$.
\end{lemma}

In almost all cases of interesting parameters the size of the incidence matrix $A_{t,k}$ is too large to find a lower bound which is close to the maximum size $B_q(t,k,n)$. Only for the parameters $P_2(2,3,6)$ a $q$-packing design with size $s=77$  was constructed from the ordinary incidence matrix $A_{2,3}$ \cite{KK08a}.

In order to reduce the matrix of coefficients $A_{t,k}$ of the integer linear optimization problem such that finding a reasonable solution  comes into feasible regions we use a \emph{group of automorphisms}. This approach was first proposed by Kramer and Mesner~\cite{KM76} who constructed ordinary combinatorial $t$-designs over finite sets.

Each subgroup $G$ of the general linear group $\gl{n}{q}$ defines a group action on the set of $k$-subspaces of $\F_q^n$ by
\[
G\times \subspaces{n}{k}{q}\to\subspaces{n}{k}{q},\, (g,K)\mapsto gK:=\{gv\mid v\in K\}, 
\]
where the $G$-orbit of $K$ will be denoted by 
\[
G(K):=\{gK\mid g\in G\}.
\]
In addition the number of all $G$-orbits on the set of $k$-subspaces will be abbreviated by $\gauss{n}{k}{q}^G$.

The $G$-incidence matrix $A_{t,k}^G=(a_{TK}^G)$ is now defined between the $G$-orbits of $t$-subspaces and the $G$-orbits of $k$-subspaces. The entry $a_{TK}^G$ corresponding to the $G$-orbit of $T$ and the $G$-orbit of $K$ is given by the number of orbit elements of $G(K)$ which contain the representative $T$:
\[
a_{TK}^G:=|\{K'\in G(K)\mid T\subseteq K'\}|.
\]

Analogous to the ordinary optimization problem we can also transform the construction task of $q$-packing designs with given cardinality admitting a group of automorphism into an integer linear programming by replacing the ordinary incidence matrix by the $G$-incidence matrix:

\begin{theorem}[Kramer-Mesner~\cite{KK08a,KM76}]
The set of $P_q(t,k,n)$ $q$-packing designs with cardinality $s$ admitting $G\le\gl{n}{q}$ as group of automorphisms bijectively corresponds to the set of 0-1-vectors $x$ such that $A_{t,k}^Gx$ is a 0-1-vector with orbit sum $\sum_K |G(K)|x_K = s$. If $x$ denotes such a 0-1-vector the corresponding $q$-packing design $\B$ is given by
\[
\B=\bigcup_{K:\, x_K=1}G(K).
\]
\end{theorem}

Note that a $P_q(t,k,n)$ $q$-packing design admitting a group of automorphisms $G\le \gl{n}{q}$ is a collection of full orbits of $G$ on $\subspaces{n}{k}{q}$. 

\section{Refinement by Subgroups}

The major drawback of the Kramer-Mesner-method is evident: For interesting parameters the reduction of the incidence matrix must be sufficiently large such that solving the optimization problem becomes feasible. But a large reduction also means that the desired $q$-packing design must admit a large group of automorphisms. In that case no $q$-packing designs with a small or even the trivial group of automorphisms can be considered.

In order to weaken the property that the complete $q$-packing design must possess the whole group of automorphisms we modify the Kramer-Mesner-approach and require that just a part of the $q$-packing design has to admit the group $G$ of automorphisms. The remaining parts of the $q$-packing designs should arise from subgroups of $G$---also from the trivial one group. 

For subgroups $H\le G\le \gl{n}{q}$ it is obvious that orbits of $G$ on $\subspaces{n}{k}{q}$ split into orbits of $H$ on $\subspaces{n}{k}{q}$. Furthermore, if we sum up the columns of the $H$-incidence matrix $A_{t,k}^H$ which correspond to the same $G$-orbit we obtain an intermediate matrix $A'$ with some identical rows: All rows in $A'$ corresponding to $H$-orbits on $t$-subspaces that fuse to the same $G$-orbit $G(T)$ on $t$-subspaces are identical due to the incidence preserving property
\[
T\subseteq K\iff gT\subseteq gK.
\]
Moreover, these rows are all identical to the row of $A_{t,k}^G$ corresponding to the orbit $G(T)$. Finally, if we erase all these identical rows in $A'$ corresponding to the same $G$-orbit and just leave exactly one row per $G$-orbit we get the $G$-incidence matrix $A_{t,k}^G$.

\begin{example}\em
We take the groups $H \le G \le \gl{4}{2}$ where $G$ is a group of order $6$ generated by
\[
\begin{bmatrix}
0&1&1&0\\
1&1&1&1\\
0&0&0&1\\
0&0&1&0
\end{bmatrix}
\]
and $H$ is the unique subgroup of $G$ of order $3$. We consider the incidence matrices $A_{1,2}^H$ and $A_{1,2}^G$ of both groups between the orbits on $1$- and $2$-subspaces. 

First we take a close look at the incidence matrix $A_{1,2}^H$. It has $7$ rows and $13$ columns which correspond to the $7$ orbits of $H$ on the $1$-subspaces and to the $13$ orbits of $H$ on the $2$-subspaces. The additional row below the doubled horizontal line indicates the sizes of the corresponding orbits. The vertical lines give a partition of all $H$-orbits on $2$-subspaces such that all orbits within the same part fuse to a $G$-orbit:
\[
A_{1,2}^H=
\left[
\begin{array}{cc|cc|c|c|cc|cc|c|c|c}
3&0&3&0&1&0&0&0&0&0&0&0&0\\
0&3&0&3&1&0&0&0&0&0&0&0&0\\
1&0&0&1&0&1&2&0&1&1&0&0&0\\
0&1&1&0&0&1&0&2&1&1&0&0&0\\
1&1&0&0&0&0&0&0&1&1&1&2&0\\
0&0&1&1&0&0&1&1&0&0&1&1&1\\
0&0&0&0&1&3&0&0&0&0&3&0&0\\
\hline
\hline
3&3&3&3&1&3&3&3&3&3&3&3&1\\
\end{array}\right]
\]
Adding up all columns within the same part yields the intermediate matrix $A'$. The horizontal lines form a partition of $H$-orbits on $1$-subspaces into parts of $H$-orbits which belong the same $G$-orbit:
\[
A^\prime=
\left[
\begin{array}{ccccccccc}
3&3&1&0&0&0&0&0&0\\
3&3&1&0&0&0&0&0&0\\
\hline
1&1&0&1&2&2&0&0&0\\
1&1&0&1&2&2&0&0&0\\
\hline
2&0&0&0&0&2&1&2&0\\
\hline
0&2&0&0&2&0&1&1&1\\
\hline
0&0&1&3&0&0&3&0&0\\
\hline
\hline
6&6&1&3&6&6&3&3&1
\end{array}\right]
\]
Deleting duplicates of the rows for the same $G$-orbit finally yields the incidence matrix $A_{1,2}^G$:
\[
A_{1,2}^G=
\left[
\begin{array}{ccccccccc}
3&3&1&0&0&0&0&0&0\\
1&1&0&1&2&2&0&0&0\\
2&0&0&0&0&2&1&2&0\\
0&2&0&0&2&0&1&1&1\\
0&0&1&3&0&0&3&0&0\\
\hline
\hline
6&6&1&3&6&6&3&3&1\\
\end{array}\right]
\]
\end{example}

The property that the incidence matrix $A_{t,k}^G$ arises from the incidence matrix $A_{t,k}^H$ by adding up columns and erasing identical rows enables us to translate a solution of the Kramer-Mesner optimization problem with prescribed group $G$ to a corresponding solution with prescribed group $H$.

\begin{lemma}
Let $H \le G \le \gl{n}{q}$ be a chain of subgroups of the general linear group, let $\{L_i\mid 1\le i \le \gauss{n}{k}{q}^H\}$ and $\{K_i\mid 1\le i \le \gauss{n}{k}{q}^G\}$, respectively, denote representatives of the orbits of $H$ and of $G$, respectively, on the set of $k$-subspaces of $\F_q^n$. Furthermore, let 
\[
x\in\{0,1\}^{\gauss{n}{k}{q}^G}\quad\text{such that}\quad
A_{t,k}^G x\in\{0,1\}^{\gauss{n}{t}{q}^G}.
\]
Since every $H$-orbit is contained in exactly one $G$-orbit, the
vector
\[
y\in\{0,1\}^{\gauss{n}{k}{q}^H},\quad y_i:=x_j\text{ if }H(L_i)\subseteq G(K_j)
\]
is well-defined and satisfies
\[
A_{t,k}^H y\in\{0,1\}^{\gauss{n}{t}{q}^H}
\qquad\text{with}\qquad
\sum_{i=1}^{\gauss{n}{k}{q}^H}|H(L_i)|y_i=\sum_{i=1}^{\gauss{n}{k}{q}^G}|G(K_i)|x_i.
\]
In addition, the $q$-packing designs defined by $x$ and $y$ are the same.
\end{lemma}

This lemma leads us to the construction approach for $q$-packing designs, which we call \emph{zoomed-Kramer-Mesner-method} in the following. The basic idea is straightforward: 
\begin{enumerate}
\item[(Z1)] 
Find a solution $x$ for $A^G_{t,k}$. 

\item[(Z2)]
Transform the found solution $x$ into a solution $y$ for $A^H_{t,k}$. 

\item[(Z3)]
Extend the solution $y$ for the matrix $A^H_{t,k}$ to $z$.
\end{enumerate}

\begin{example}[cont.]\em
Considering the last example a maximal $P_2(1,2,4)$ $q$-packing design admitting $G$ as group of automorphisms has cardinality $s=2$ and it is defined by the vector $x = 001000001$:
\[
A_{1,2}^G=
\left[
\begin{array}{ccccccccc}
3&3&\bf 1&0&0&0&0&0&\bf0\\
1&1&\bf0&1&2&2&0&0&\bf0\\
2&0&\bf0&0&0&2&1&2&\bf0\\
0&2&\bf0&0&2&0&1&1&\bf1\\
0&0&\bf1&3&0&0&3&0&\bf0\\
\hline
\hline
6&6&\bf1&3&6&6&3&3&\bf1\\
\end{array}\right]
\]
Transforming $x$ into a solution for the subgroup $H\le G$ yields the vector $y=0000100000001$:
\[
A_{1,2}^H=
\left[
\begin{array}{cc|cc|c|c|cc|cc|c|c|c}
3&0&3&0&\bf1&0&0&0&0&0&0&0&\bf0\\
0&3&0&3&\bf1&0&0&0&0&0&0&0&\bf0\\
1&0&0&1&\bf0&1&2&0&1&1&0&0&\bf0\\
0&1&1&0&\bf0&1&0&2&1&1&0&0&\bf0\\
1&1&0&0&\bf0&0&0&0&1&1&1&2&\bf0\\
0&0&1&1&\bf0&0&1&1&0&0&1&1&\bf1\\
0&0&0&0&\bf1&3&0&0&0&0&3&0&\bf0\\
\hline
\hline
3&3&3&3&\bf1&3&3&3&3&3&3&3&\bf1\\
\end{array}\right]
\]
A possible extension of $y$ is $z=0000100001001$ which finally defines a $q$-packing design of size $s=5$.
\end{example}

This zoomed-Kramer-Mesner-approach can be refined using the following strategies:

\subsection*{Traversing subgroup chains}

A possible strategy is to take a chain of subgroups of the general linear group $G=G_1\ge G_2\ge \ldots\ge G_\alpha$ and step with the search through all incidence matrices. With large groups $G_i$ at the beginning we have an incidence matrix $A_{t,k}^{G_i}$ of moderate size such that the optimization problems remains feasible. The refinement and extension of solutions along a chain of subgroups can be seen as a kind of \lq\lq zooming\rq\rq{} into the incidence matrix $A_{t,k}^G$ where the resolution increases from a subgroup to a smaller subgroup.

\subsection*{Cutting off the Kramer-Mesner-matrix}

In fact when we try to find an admissible solution on the Kramer-Mesner-matrix $A_{t,k}^G$ all columns with entries greater than one can be erased before starting the search. The issue is now: if we find a maximal solution $x$ for the group $G$ and if we want to extend this solution for a subgroup $H$, are there any $H$-orbits that can directly be omitted for the continued search? To answer this question we consider the following observation.

\begin{lemma}\label{lem:fitin}
Let $\B$ be a $P_q(t,k,n)$ $q$-packing design admitting $G\le \gl{n}{q}$ as a group of automorphisms and let $K$ be a $k$-subspace of $\F_q^n$ such that $\dim(K\cap K')<t$ for all $K'\in \B$. Then, we have $\dim(gK\cap K')<t$ for all $K'\in \B$ and $g\in G$.
\end{lemma}

We call a $G$-orbit on $k$-subspaces \emph{admissible} if the corresponding column in the Kramer-Mesner-matrix $A_{t,k}^G$ has only entries less than or equal to one. It is clear that a $q$-packing design $\B$ admitting $G$ as group of automorphisms only consists of admissible orbits. 

With Lemma~\ref{lem:alternative} and \ref{lem:fitin} we can directly follow that if $\B$ can be extended by any element of any admissible orbit $G(K)$ we can finally extend $\B$ by the \emph{complete} orbit $G(K)$. Hence, if $\B$ is maximal (with respect to extension by $G$-orbits) no $k$-subspace that is contained in any admissible $G$-orbit can be added to $\B$. The extension of a maximal $\B$ by orbits of subgroups $H\le G$ can only be done by $H$-orbits that are contained in \emph{non-admissible} $G$-orbits. Therefore, all $H$-orbits which are parts of admissible $G$-orbits can be omitted for further consideration.

\subsection*{Local modifications}

If we transform a solution $x$ admitting a group $G$ to a solution $y$ admitting a subgroup $H$, it can be advantageous first to modify the solution $y$ by deleting or exchanging already selected columns in $y$ and to try to extend the altered vector $y'$ instead of the original $y$. Using these modifications $H$-orbits which are contained in admissible $G$-orbits can be reconsidered for further extension.

\section{Metaheuristic Optimization}

Heuristic solvers are approximative methods to tackle large optimization problems with the goal to produce satisfiable solutions---not necessarily the best possible solution. Heuristics which can be easily adapted to different problems are called metaheuristics. Our metaheuristic solver is based on a well known breadth-first search method named \emph{Beam-search}~\cite{Bis87}. 

The algorithm takes parameters $\alpha$, $\beta$, and an objective function $f$. The integer $\alpha$ is the predetermined number of solutions (which will be called $\alpha$-solutions) at each construction step to be examined in parallel, also called the beam width. The integer $\beta$ is the number of extensions to be considered for each of the $\alpha$-solutions. In our case one $\alpha$-solution represents a 0-1-solution vector $x$ such that $A_{t,k}^G x$ yields a 0-1-vector. Each iteration of our search method works as follows: 

\begin{enumerate}
\item[(B1)] 
During the first iteration step we choose randomly $\alpha$ different columns of $A_{t,k}^G$ (which are called states).

\item[(B2)]
For every of those $\alpha$ states we eliminate all columns of $A_{t,k}^G$ which have at least one nonzero entry in the same row as the chosen state. We keep the amount of remaining columns for each state as $f$. 

\item[(B3)]
In all other iteration steps we choose maximal $\alpha\beta$ extensions uniformly at random and restrict us to the $\alpha$ solutions with highest $f$-values---these solutions are the new $\alpha$-solutions. For any non-feasible $\alpha$-solution, we set $\alpha-1$. The iteration ends when $\alpha=0$ is reached. The whole search ends when some predefined condition is satisfied. 
\end{enumerate}

To handle the elimination step we use the following data structures: We only hold the row indices with entries greater than zero for each column of the incidence matrix $A_{t,k}^G$ in a binary tree such that we get a vector $\Delta$ of binary trees. Every $\alpha$-solution is associated with a vector containing the chosen states and a vector holding the remaining possible states. After choosing a state of the $\alpha\beta$ extensions we eliminate all remaining states covering any row indices of the chosen state in $\Delta$. We finally add the current state to the already chosen states of the corresponding $\alpha$-solution.

\section{Results}

Using the proposed methods (Kramer-Mesner-approach, zoomed-Kramer-Mesner-approach, Beam-search) we can improve some binary $q$-packing designs. A summary of all new results is given by Table~\ref{tab:new}.

\begin{table}[!htbp]
\centering
\caption{New lower bounds on $B_2(2,3,n)$}\label{tab:new}
\begin{tabular}{lrrlrcc}
\toprule
$n$ & Old  & New  & Group& Order &Zoom & Local Mod.\\
\midrule
$7$ & $304$ & $329$ & cyclic & $15$ & \checkmark & \checkmark \\
$8$ & $1275$ & $1312$ & cyclic & $217$ & \checkmark & \\
$11$& $79833$ & $92411$ &$\ns{10}{2}\times1$& $22517$ &&\\
$12$& $315315$ & $385515$ & $\ns{12}{2}$ & $49140$&&\\
$14$& $4177665$ &$5996178$ &$\ns{14}{2}$ & $229362$&&\\
\bottomrule
\end{tabular}
\end{table}

In order to get a convenient representation of subspaces we display a $3$-subspace generated by the three column vectors
\[
\begin{bmatrix}
x_0 & y_0 & z_0\\
x_1 & y_1 & z_1\\
\vdots &\vdots&\vdots\\
x_{n-1} & y_{n-1} & z_{n-1}\\
\end{bmatrix}
\]
as a triple of integers 
\[
[X,Y,Z]=\left[\sum_{i=0}^{n-1}x_i2^i,\sum_{i=0}^{n-1}y_i2^i, \sum_{i=0}^{n-1}z_i2^i\right].
\]

Furthermore, by $\ns{n}{q}$ we denote the \emph{normalizer of a Singer cycle}, which is generated by a Singer cycle the Frobenius automorphism~\cite{BEO+12}. The order of $\ns{n}{q}$ is $n(q^n-1)$.

\subsection*{$B_2(2,3,7)\ge 329$}

We prescribe a group $G$ of order $15$ generated by
\[\footnotesize
\arraycolsep=3pt
\left[\begin{array}{ccccccc}
0& 1& 0& 1& 0& 0& 0\\ 
1& 0& 0& 0& 1& 1& 0\\ 
1& 1& 0& 0& 1& 1& 1\\
1& 1& 0& 0& 0& 0& 0\\
1& 1& 1& 0& 0& 1& 1\\
1& 0& 0& 1& 1& 1& 0\\
0& 1& 0& 1& 0& 1& 1
\end{array}\right].
\]
Applying the zoomed version of Kramer-Mesner-method with the trivial subgroup $1=H\le G$ and breaking further symmetries by exchanging selected columns we are able to find a $P_2(2,3,7)$ $q$-packing design with size $329$. The set of all of these $3$-subspaces is given in Table~\ref{tab:2_3_7}.

\begin{table}[!htbp]
\caption{$P_2(2,3,7)$ $q$-packing design with size $329$}\label{tab:2_3_7}
\centering
{\scriptsize
\tabcolsep=1pt
\begin{tabular}{lllllll}
\toprule
$[7,32,64]$,&  $[15,16,64]$,&  $[18,8,64]$,&  $[3,40,64]$,&  $[54,24,64]$,&  $[1,56,64]$,&  $[17,4,64]$,\\
$[42,36,64]$,&  $[41,52,64]$,&  $[33,12,64]$,&  $[9,44,64]$,&  $[51,28,64]$,&  $[2,60,64]$,&  $[25,34,64]$,\\
$[5,50,64]$,&  $[49,58,64]$,&  $[53,22,64]$,&  $[13,30,64]$,&  $[14,80,32]$,&  $[83,72,32]$,&  $[85,24,32]$,\\
$[69,88,32]$,&  $[18,4,32]$,&  $[90,68,32]$,&  $[17,20,32]$,&  $[67,84,32]$,&  $[91,12,32]$,&  $[89,76,32]$,\\
$[65,28,32]$,&  $[26,92,32]$,&  $[9,66,32]$,&  $[25,10,32]$,&  $[73,6,32]$,&  $[81,78,32]$,&  $[65,16,96]$,\\
$[76,80,96]$,&  $[4,8,96]$,&  $[6,72,96]$,&  $[3,24,96]$,&  $[9,68,96]$,&  $[26,20,96]$,&  $[18,84,96]$,\\
$[17,2,96]$,&  $[13,66,96]$,&  $[89,82,96]$,&  $[29,10,96]$,&  $[85,74,96]$,&  $[1,90,96]$,&  $[73,22,96]$,\\
$[5,86,96]$,&  $[21,94,96]$,&  $[68,8,16]$,&  $[3,72,16]$,&  $[66,40,16]$,&  $[70,104,16]$,&  $[105,4,16]$,\\
$[1,100,16]$,&  $[9,12,16]$,&  $[42,44,16]$,&  $[11,108,16]$,&  $[69,2,16]$,&  $[97,34,16]$,&  $[45,10,16]$,\\
$[41,74,16]$,&  $[13,102,16]$,&  $[37,14,16]$,&  $[73,78,16]$,&  $[77,110,16]$,&  $[5,8,80]$,&  $[101,72,80]$,\\
$[36,40,80]$,&  $[66,68,80]$,&  $[67,100,80]$,&  $[10,108,80]$,&  $[73,2,80]$,&  $[37,34,80]$,&  $[1,98,80]$,\\
$[65,74,80]$,&  $[105,106,80]$,&  $[33,70,80]$,&  $[109,78,80]$,&  $[41,110,80]$,&  $[33,72,48]$,&  $[1,40,48]$,\\
$[71,104,48]$,&  $[9,4,48]$,&  $[75,68,48]$,&  $[10,36,48]$,&  $[3,100,48]$,&  $[67,12,48]$,&  $[65,108,48]$,\\
$[37,2,48]$,&  $[73,42,48]$,&  $[109,106,48]$,&  $[77,70,48]$,&  $[97,14,48]$,&  $[101,78,48]$,&  $[5,110,48]$,\\
$[71,8,112]$,&  $[67,72,112]$,&  $[97,40,112]$,&  $[2,104,112]$,&  $[35,68,112]$,&  $[107,100,112]$,&  $[98,12,112]$,\\
$[105,76,112]$,&  $[41,44,112]$,&  $[34,108,112]$,&  $[45,66,112]$,&  $[9,10,112]$,&  $[101,74,112]$,&  $[13,42,112]$,\\
$[1,6,112]$,&  $[109,38,112]$,&  $[51,36,8]$,&  $[81,52,8]$,&  $[65,116,8]$,&  $[97,2,8]$,&  $[33,50,8]$,\\
$[49,70,8]$,&  $[113,38,8]$,&  $[1,102,8]$,&  $[21,22,8]$,&  $[17,54,8]$,&  $[37,118,8]$,&  $[34,4,72]$,\\
$[97,36,72]$,&  $[99,20,72]$,&  $[113,84,72]$,&  $[5,66,72]$,&  $[117,98,72]$,&  $[85,50,72]$,&  $[81,114,72]$,\\
$[49,22,72]$,&  $[1,86,72]$,&  $[114,4,40]$,&  $[115,68,40]$,&  $[35,100,40]$,&  $[19,20,40]$,&  $[33,84,40]$,\\
$[5,18,40]$,&  $[101,50,40]$,&  $[37,70,40]$,&  $[21,38,40]$,&  $[69,22,40]$,&  $[113,86,40]$,&  $[65,54,40]$,\\
$[113,68,104]$,&  $[81,20,104]$,&  $[34,52,104]$,&  $[101,98,104]$,&  $[33,18,104]$,&  $[117,82,104]$,&  $[17,50,104]$,\\
$[21,6,104]$,&  $[65,38,104]$,&  $[5,102,104]$,&  $[97,54,104]$,&  $[19,4,24]$,&  $[98,36,24]$,&  $[97,100,24]$,\\
$[50,20,24]$,&  $[115,52,24]$,&  $[35,116,24]$,&  $[113,66,24]$,&  $[69,34,24]$,&  $[65,18,24]$,&  $[37,82,24]$,\\
$[33,6,24]$,&  $[81,102,24]$,&  $[17,22,24]$,&  $[21,86,24]$,&  $[99,68,88]$,&  $[113,100,88]$,&  $[2,20,88]$,\\
$[83,84,88]$,&  $[50,52,88]$,&  $[82,116,88]$,&  $[37,66,88]$,&  $[53,98,88]$,&  $[49,18,88]$,&  $[1,114,88]$,\\
$[117,70,88]$,&  $[65,86,88]$,&  $[99,4,56]$,&  $[50,68,56]$,&  $[35,36,56]$,&  $[81,100,56]$,&  $[2,84,56]$,\\
$[113,52,56]$,&  $[83,116,56]$,&  $[21,98,56]$,&  $[37,18,56]$,&  $[5,6,56]$,&  $[17,70,56]$,&  $[49,38,56]$,\\
$[33,102,56]$,&  $[101,22,56]$,&  $[83,4,120]$,&  $[49,68,120]$,&  $[98,100,120]$,&  $[51,20,120]$,&  $[50,84,120]$,\\
$[66,52,120]$,&  $[5,2,120]$,&  $[53,34,120]$,&  $[81,18,120]$,&  $[113,82,120]$,&  $[21,114,120]$,&  $[85,38,120]$,\\
$[97,22,120]$,&  $[37,54,120]$,&  $[49,2,4]$,&  $[97,66,4]$,&  $[113,10,4]$,&  $[1,74,4]$,&  $[41,42,4]$,\\
$[73,58,4]$,&  $[57,34,68]$,&  $[97,98,68]$,&  $[65,82,68]$,&  $[33,42,68]$,&  $[1,106,68]$,&  $[41,26,68]$,\\
$[121,58,68]$,&  $[81,2,36]$,&  $[49,66,36]$,&  $[65,34,36]$,&  $[17,82,36]$,&  $[121,114,36]$,&  $[73,74,36]$,\\
$[113,26,36]$,&  $[33,90,36]$,&  $[89,2,100]$,&  $[57,18,100]$,&  $[33,82,100]$,&  $[49,10,100]$,&  $[25,74,100]$,\\
$[65,58,100]$,&  $[41,34,20]$,&  $[65,98,20]$,&  $[1,82,20]$,&  $[113,114,20]$,&  $[25,106,20]$,&  $[81,34,84]$,\\
$[121,82,84]$,&  $[17,10,84]$,&  $[65,42,84]$,&  $[49,106,84]$,&  $[97,26,84]$,&  $[73,90,84]$,&  $[89,122,84]$,\\
$[97,82,52]$,&  $[1,10,52]$,&  $[33,26,52]$,&  $[65,90,52]$,&  $[25,122,52]$,&  $[33,34,116]$,&  $[57,98,116]$,\\
$[25,50,116]$,&  $[89,74,116]$,&  $[113,106,116]$,&  $[81,90,116]$,&  $[9,58,116]$,&  $[121,66,12]$,&  $[17,34,12]$,\\
$[81,82,12]$,&  $[57,114,12]$,&  $[41,10,12]$,&  $[65,106,12]$,&  $[1,26,12]$,&  $[49,90,12]$,&  $[89,58,12]$,\\
$[105,122,12]$,&  $[73,98,76]$,&  $[113,18,76]$,&  $[41,82,76]$,&  $[49,50,76]$,&  $[65,114,76]$,&  $[81,74,76]$,\\
$[17,42,76]$,&  $[9,26,76]$,&  $[57,122,76]$,&  $[65,2,44]$,&  $[17,66,44]$,&  $[1,34,44]$,&  $[105,98,44]$,\\
$[121,106,44]$,&  $[49,26,44]$,&  $[97,122,44]$,&  $[57,66,108]$,&  $[49,98,108]$,&  $[1,18,108]$,&  $[17,114,108]$,\\
$[33,74,108]$,&  $[105,42,108]$,&  $[25,26,108]$,&  $[81,122,108]$,&  $[1,66,28]$,&  $[121,18,28]$,&  $[105,50,28]$,\\
$[81,42,28]$,&  $[33,106,28]$,&  $[89,90,28]$,&  $[17,58,28]$,&  $[41,122,28]$,&  $[57,2,92]$,&  $[25,66,92]$,\\
$[33,98,92]$,&  $[17,18,92]$,&  $[97,114,92]$,&  $[49,42,92]$,&  $[73,106,92]$,&  $[9,90,92]$,&  $[1,122,92]$,\\
$[9,18,60]$,&  $[33,114,60]$,&  $[121,74,60]$,&  $[17,106,60]$,&  $[57,26,60]$,&  $[49,82,124]$,&  $[41,114,124]$,\\
$[113,74,124]$,&  $[89,106,124]$,&  $[25,90,124]$,&  $[17,122,124]$, & $[2,68,72]$,& $[69,6,8]$,&$[89,26,4]$,\\
\bottomrule
\end{tabular}}
\end{table}

\subsection*{$B_2(2,3,8)\ge 1312$}

We use a subgroup chain $H\le G$ where $G$ is a group of order $217$ generated by
\[\footnotesize
\arraycolsep=3pt
\left[\begin{array}{cccccccc}
0& 0& 1& 1& 0& 0& 0& 0\\ 
1& 0& 0& 0& 1& 1& 0& 0\\ 
0& 0& 0& 0& 0& 0& 1& 1\\ 
0& 1& 1& 1& 1& 0& 1& 1\\ 
0& 1& 0& 0& 1& 1& 0& 1\\ 
1& 1& 0& 1& 0& 0& 0& 0\\ 
0& 0& 1& 0& 0& 1& 1& 0\\ 
0& 0& 1& 0& 1& 0& 0& 1 
\end{array}\right]
\]
and $H$ is the subgroup of order $7$. The zoomed-Kramer-Mesner-approach without any local modifications yields a $q$-packing design with size $1312$. The $196$ representatives of $H$-orbits are listed in Table~\ref{tab:2_3_8}.

\begin{table}[!htbp]
\caption{$P_2(2,3,8)$ $q$-packing design with size $1312$}\label{tab:2_3_8}
\centering
{\scriptsize
\tabcolsep=1pt
\begin{tabular}{llllll}
\toprule
$[50,64,128]$, &  $[90,32,128]$, &  $[101,16,128]$, &  $[65,80,128]$, &  $[1,72,128]$, &  $[21,40,128]$,\\
$[118,24,128]$, &  $[6,88,128]$, &  $[23,56,128]$, &  $[20,120,128]$, &  $[82,4,128]$, &  $[25,68,128]$,\\
$[2,36,128]$, &  $[19,100,128]$, &  $[113,84,128]$, &  $[67,116,128]$, &  $[57,12,128]$, &  $[98,76,128]$,\\
$[43,44,128]$, &  $[75,28,128]$, &  $[35,92,128]$, &  $[115,124,128]$, &  $[125,34,128]$, &  $[13,18,128]$,\\
$[9,10,128]$, &  $[45,74,128]$, &  $[49,42,128]$, &  $[85,106,128]$, &  $[33,26,128]$, &  $[77,70,128]$,\\
$[89,22,128]$, &  $[5,54,128]$, &  $[97,14,128]$, &  $[109,62,128]$, &  $[8,32,64]$, &  $[132,160,64]$,\\
$[140,16,64]$, &  $[169,176,64]$, &  $[149,136,64]$, &  $[54,168,64]$, &  $[22,184,64]$, &  $[33,4,64]$,\\
$[170,164,64]$, &  $[35,20,64]$, &  $[146,148,64]$, &  $[145,52,64]$, &  $[43,180,64]$, &  $[3,12,64]$,\\
$[139,44,64]$, &  $[19,172,64]$, &  $[2,60,64]$, &  $[163,188,64]$, &  $[185,130,64]$, &  $[173,42,64]$,\\
$[41,26,64]$, &  $[45,134,64]$, &  $[49,38,64]$, &  $[17,166,64]$, &  $[53,142,64]$, &  $[161,46,64]$,\\
$[57,30,64]$, &  $[129,32,192]$, &  $[29,160,192]$, &  $[132,16,192]$, &  $[2,144,192]$, &  $[137,48,192]$,\\
$[142,176,192]$, &  $[180,8,192]$, &  $[33,136,192]$, &  $[135,40,192]$, &  $[34,168,192]$, &  $[52,152,192]$,\\
$[182,184,192]$, &  $[9,4,192]$, &  $[154,36,192]$, &  $[42,12,192]$, &  $[145,28,192]$, &  $[25,10,192]$,\\
$[165,170,192]$, &  $[17,6,192]$, &  $[177,134,192]$, &  $[45,166,192]$, &  $[57,22,192]$, &  $[49,54,192]$,\\
$[173,174,192]$, &  $[21,30,192]$, &  $[37,158,192]$, &  $[132,144,32]$, &  $[67,80,32]$, &  $[76,208,32]$,\\
$[131,88,32]$, &  $[201,68,32]$, &  $[146,196,32]$, &  $[154,84,32]$, &  $[153,212,32]$, &  $[1,204,32]$,\\
$[147,28,32]$, &  $[9,92,32]$, &  $[11,220,32]$, &  $[217,18,32]$, &  $[137,210,32]$, &  $[69,70,32]$,\\
$[145,142,32]$, &  $[73,94,32]$, &  $[13,222,32]$, &  $[66,16,160]$, &  $[136,144,160]$, &  $[131,80,160]$,\\
$[7,8,160]$, &  $[81,72,160]$, &  $[219,20,160]$, &  $[193,84,160]$, &  $[91,12,160]$, &  $[27,140,160]$,\\
$[138,76,160]$, &  $[210,204,160]$, &  $[209,28,160]$, &  $[146,92,160]$, &  $[129,202,160]$, &  $[137,6,160]$,\\
$[213,22,160]$, &  $[77,214,160]$, &  $[153,142,160]$, &  $[13,158,160]$, &  $[17,94,160]$, &  $[65,222,160]$,\\
$[143,16,96]$, &  $[139,208,96]$, &  $[132,136,96]$, &  $[148,200,96]$, &  $[81,4,96]$, &  $[82,68,96]$,\\
$[147,196,96]$, &  $[1,20,96]$, &  $[130,84,96]$, &  $[195,140,96]$, &  $[193,18,96]$, &  $[129,210,96]$,\\
$[145,138,96]$, &  $[205,6,96]$, &  $[65,150,96]$, &  $[141,86,96]$, &  $[213,78,96]$, &  $[137,206,96]$,\\
$[9,158,96]$, &  $[149,222,96]$, &  $[204,144,224]$, &  $[79,80,224]$, &  $[67,208,224]$, &  $[198,136,224]$,\\
$[68,88,224]$, &  $[18,20,224]$, &  $[91,212,224]$, &  $[25,12,224]$, &  $[90,220,224]$, &  $[217,138,224]$,\\
$[193,74,224]$, &  $[153,14,224]$, &  $[85,30,224]$, &  $[201,158,224]$, &  $[33,8,16]$, &  $[73,164,16]$,\\
$[202,228,16]$, &  $[194,12,16]$, &  $[201,170,16]$, &  $[1,38,16]$, &  $[97,36,144]$, &  $[33,76,144]$,\\
$[137,226,144]$, &  $[1,8,80]$, &  $[73,172,80]$, &  $[138,108,80]$, &  $[173,38,80]$, &  $[70,136,208]$,\\
$[10,164,208]$, &  $[137,204,208]$, &  $[201,194,208]$, &  $[99,104,48]$, &  $[35,4,48]$, &  $[75,68,48]$,\\
$[161,228,48]$, &  $[227,168,176]$, &  $[69,134,176]$, &  $[97,40,112]$, &  $[202,196,112]$, &  $[42,100,240]$,\\
$[229,166,240]$, &  $[18,228,168]$, &  $[113,148,232]$, &  $[129,86,232]$, &  $[67,100,24]$, &  $[229,50,24]$,\\
$[69,226,88]$, &  $[105,194,100]$, &  $[153,194,140]$, &  $[113,218,140]$\\
\bottomrule
\end{tabular}}
\end{table}

\subsection*{$B_2(2,3,11)\ge 92411$}

We prescribed the group $G =\ns{10}{2}\times{1}\le\gl{11}{2}$ of order $10230$ generated by the Singer cycle and the Frobenius automorphism of dimension $10$ extended by a further dimension:
\[
\footnotesize
\arraycolsep=3pt
\left[\begin{array}{cccccccccc|c}
0& 0& 0& 0& 0& 0& 0& 0& 0& 1& 0\\ 
1& 0& 0& 0& 0& 0& 0& 0& 0& 0& 0\\ 
0& 1& 0& 0& 0& 0& 0& 0& 0& 0& 0\\
0& 0& 1& 0& 0& 0& 0& 0& 0& 1& 0\\ 
0& 0& 0& 1& 0& 0& 0& 0& 0& 0& 0\\ 
0& 0& 0& 0& 1& 0& 0& 0& 0& 0& 0\\ 
0& 0& 0& 0& 0& 1& 0& 0& 0& 0& 0\\ 
0& 0& 0& 0& 0& 0& 1& 0& 0& 0& 0\\ 
0& 0& 0& 0& 0& 0& 0& 1& 0& 0& 0\\ 
0& 0& 0& 0& 0& 0& 0& 0& 1& 0& 0\\ 
\hline
0& 0& 0& 0& 0& 0& 0& 0& 0& 0& 1 
\end{array}\right]\,,\,
\left[\begin{array}{cccccccccc|c}
1&0&0&0&0&1&0&0&0&0&0\\
0&0&0&0&0&0&0&0&0&1&0\\
0&1&0&0&0&0&1&0&0&0&0\\
0&0&0&0&0&1&0&0&0&0&0\\
0&0&1&0&0&0&0&1&0&1&0\\
0&0&0&0&0&0&1&0&0&0&0\\
0&0&0&1&0&0&0&0&1&0&0\\
0&0&0&0&0&0&0&1&0&0&0\\
0&0&0&0&1&0&0&0&0&1&0\\
0&0&0&0&0&0&0&0&1&0&0\\
\hline
0&0&0&0&0&0&0&0&0&0&1
\end{array}\right].
\]
We are able to find many solutions of a $P_2(2,3,11)$ $q$-packing design consisting of $11$ orbits with cardinality is $92411$. We used the \lq\lq plain\rq\rq{} Kramer-Mesner-approach without refinement by subgroups. The representatives of $G$-orbits are given in Table~\ref{tab:2_3_11}.

\begin{table}[!htbp]
\caption{$P_2(2,3,11)$ $q$-packing design with size $92411$}\label{tab:2_3_11}
\centering
{\scriptsize
\begin{tabular}{llll}
\toprule
$[454,512,1024]$,& $[100,256,512]$,& $[239,1280,512]$,& $[260,1408,512]$,\\
$[1460,64,512]$,& $[278,1088,512]$,& $[23,1216,512]$,& $[1471,448,512]$,\\
$[1324,1040,512]$,& $[1383,1488,512]$,& $[232,1264,512]$\\
\bottomrule
\end{tabular}}
\end{table}

\subsection*{$B_2(2,3,12)\ge 385515$}

We prescribed the normalizer of the Singer cycle $\ns{12}{2}$ which has order $49140$ and which is generated by
\[
\footnotesize
\arraycolsep=3pt
\left[\begin{array}{cccccccccccc}
0& 0& 0& 0& 0& 0& 0& 0& 0& 0& 0& 1\\ 
1& 0& 0& 0& 0& 0& 0& 0& 0& 0& 0& 1\\ 
0& 1& 0& 0& 0& 0& 0& 0& 0& 0& 0& 0\\
0& 0& 1& 0& 0& 0& 0& 0& 0& 0& 0& 0\\
0& 0& 0& 1& 0& 0& 0& 0& 0& 0& 0& 1\\
0& 0& 0& 0& 1& 0& 0& 0& 0& 0& 0& 0\\
0& 0& 0& 0& 0& 1& 0& 0& 0& 0& 0& 1\\
0& 0& 0& 0& 0& 0& 1& 0& 0& 0& 0& 0\\
0& 0& 0& 0& 0& 0& 0& 1& 0& 0& 0& 0\\
0& 0& 0& 0& 0& 0& 0& 0& 1& 0& 0& 0\\
0& 0& 0& 0& 0& 0& 0& 0& 0& 1& 0& 0\\
0& 0& 0& 0& 0& 0& 0& 0& 0& 0& 1& 0
\end{array}\right]\,,\,
\left[\begin{array}{cccccccccccc}
1& 0& 0& 0& 0& 0& 1& 0& 0& 1& 1& 0\\
0& 0& 0& 0& 0& 0& 1& 0& 0& 1& 1& 0\\ 
0& 1& 0& 0& 0& 0& 0& 1& 0& 0& 1& 1\\
0& 0& 0& 0& 0& 0& 0& 1& 0& 0& 1& 1\\
0& 0& 1& 0& 0& 0& 1& 0& 1& 1& 1& 1\\
0& 0& 0& 0& 0& 0& 0& 0& 1& 0& 0& 1\\
0& 0& 0& 1& 0& 0& 1& 1& 0& 0& 0& 1\\
0& 0& 0& 0& 0& 0& 0& 0& 0& 1& 0& 0\\
0& 0& 0& 0& 1& 0& 0& 1& 1& 0& 0& 0\\
0& 0& 0& 0& 0& 0& 0& 0& 0& 0& 1& 0\\
0& 0& 0& 0& 0& 1& 0& 0& 1& 1& 0& 0\\
0& 0& 0& 0& 0& 0& 0& 0& 0& 0& 0& 1
\end{array}\right].
\]
The simple version of the Kramer-Mesner-method yields $10$ orbit representatives defining the $q$-packing design of size $385515$ which are listed in Table~\ref{tab:2_3_12}.

\begin{table}[!htbp]
\caption{$P_2(2,3,12)$ $q$-packing design with size $385515$}\label{tab:2_3_12}
\centering
{\scriptsize
\begin{tabular}{llll}
\toprule
$[829,1024,2048]$, & $[1306,1536,2048]$, & $[1272,256,2048]$, & $[1213,1280,2048]$,\\$[375,1152,2048]$, & $[309,384,2048]$, & $[1575,1408,2048]$, & $[1281,64,2048]$, \\
$[526,1184,2048]$, & $[2034,1060,2048]$ & & \\
\bottomrule
\end{tabular}}
\end{table}

\subsection*{$B_2(2,3,14)\ge 5996178$}
Again the prescription the normalizer of the Singer cycle $G$ of order $229362$ is used:
\[
\footnotesize
\arraycolsep=3pt
\left[\begin{array}{cccccccccccccc}
0& 0& 0& 0& 0& 0& 0& 0& 0& 0& 0& 0& 0& 1\\
1& 0& 0& 0& 0& 0& 0& 0& 0& 0& 0& 0& 0& 1\\ 
0& 1& 0& 0& 0& 0& 0& 0& 0& 0& 0& 0& 0& 0\\ 
0& 0& 1& 0& 0& 0& 0& 0& 0& 0& 0& 0& 0& 1\\ 
0& 0& 0& 1& 0& 0& 0& 0& 0& 0& 0& 0& 0& 0\\
0& 0& 0& 0& 1& 0& 0& 0& 0& 0& 0& 0& 0& 1\\
0& 0& 0& 0& 0& 1& 0& 0& 0& 0& 0& 0& 0& 0\\
0& 0& 0& 0& 0& 0& 1& 0& 0& 0& 0& 0& 0& 0\\
0& 0& 0& 0& 0& 0& 0& 1& 0& 0& 0& 0& 0& 0\\
0& 0& 0& 0& 0& 0& 0& 0& 1& 0& 0& 0& 0& 0\\
0& 0& 0& 0& 0& 0& 0& 0& 0& 1& 0& 0& 0& 0\\
0& 0& 0& 0& 0& 0& 0& 0& 0& 0& 1& 0& 0& 0\\
0& 0& 0& 0& 0& 0& 0& 0& 0& 0& 0& 1& 0& 0\\
0& 0& 0& 0& 0& 0& 0& 0& 0& 0& 0& 0& 1& 0 
\end{array}\right]\,,\,
\left[\begin{array}{cccccccccccccc}
1& 0& 0& 0& 0& 0& 0& 1& 0& 0& 0& 0& 0& 0\\ 
0& 0& 0& 0& 0& 0& 0& 1& 0& 0& 0& 0& 1& 1\\ 
0& 1& 0& 0& 0& 0& 0& 0& 1& 0& 0& 0& 1& 1\\
0& 0& 0& 0& 0& 0& 0& 1& 1& 0& 0& 0& 0& 1\\
0& 0& 1& 0& 0& 0& 0& 0& 0& 1& 0& 0& 1& 0\\
0& 0& 0& 0& 0& 0& 0& 1& 1& 1& 0& 0& 0& 0\\
0& 0& 0& 1& 0& 0& 0& 0& 0& 0& 1& 0& 1& 0\\
0& 0& 0& 0& 0& 0& 0& 0& 1& 1& 1& 0& 0& 0\\
0& 0& 0& 0& 1& 0& 0& 0& 0& 0& 0& 1& 0& 1\\
0& 0& 0& 0& 0& 0& 0& 0& 0& 1& 1& 1& 0& 0\\
0& 0& 0& 0& 0& 1& 0& 0& 0& 0& 0& 0& 1& 0\\
0& 0& 0& 0& 0& 0& 0& 0& 0& 0& 1& 1& 1& 0\\
0& 0& 0& 0& 0& 0& 1& 0& 0& 0& 0& 0& 0& 1\\
0& 0& 0& 0& 0& 0& 0& 0& 0& 0& 0& 1& 1& 1
\end{array}\right].
\]
The plain version of the Kramer-Mesner-method yields a $q$-packing design with cardinality $5996178$. The corresponding $28$ representatives of $G$-orbits are given in Table~\ref{tab:2_3_14}.

\begin{table}[!htbp]
\caption{$P_2(2,3,14)$ $q$-packing design with size $5996178$}\label{tab:2_3_14}
\centering
{\scriptsize
\begin{tabular}{llll}
\toprule
$[3175,4096,8192]$,& $[879,6144,8192]$,& $[4825,5120,8192]$,& $[2934,3072,8192]$,\\
$[47,7168,8192]$,& $[6233,4608,8192]$,& $[1222,3584,8192]$,& $[4097,7680,8192]$,\\
$[4681,4352,8192]$,& $[3300,6400,8192]$,& $[1113,5376,8192]$,& $[7797,7424,8192]$,\\
$[7198,768,8192]$,& $[5692,4864,8192]$,& $[6225,6912,8192]$,& $[1746,1792,8192]$,\\
$[255,5888,8192]$,& $[829,5248,8192]$,& $[7026,5760,8192]$,& $[6241,3456,8192]$,\\
$[5472,896,8192]$,& $[6954,7040,8192]$,& $[1852,576,8192]$,& $[6711,6720,8192]$,\\
$[3974,2368,8192]$,& $[5399,2880,8192]$,& $[1598,3264,8192]$,& $[7270,7920,8192]$\\
\bottomrule
\end{tabular}}
\end{table}
\section{Conclusion and Outlook}

The proposed zoomed-Kramer-Mesner-approach offers extended possibilities to integer linear programming on incidence structures compared to the ordinary Kramer-Mesner-approach since we give up the requirement that the complete object has to admit a given group of automorphisms. 

In the case of $q$-packing designs we were able to improve some lower bounds of $B_2(2,3,n)$. In lower dimensions we improved the bounds with the zoomed version of Kramer-Mesner-method using cyclic groups. For higher dimensions we used the plain version of the Kramer-Mesner-approach together with the normalizer of the Singer cycle which always yields an improvement of the lower bound.

Of course the computational limits of this approach have to be examined intensively. Frankly speaking we do not expect that the zoomed-Kramer-Mesner-methods works properly for large $n$ since the construction of the orbits of an arbitrarily chosen group $G\le \gl{n}{q}$ on the $t$- and $k$-subspaces exceeds computational limits. Values $n\le 15$ are realistic for $t< k\le n/2$.

Nevertheless, the zoomed version of the Kramer-Mesner-approach can also be applied to similar construction problems of incidence structures. Several articles describe the successful construction of linear block codes~\cite{BKW05,Bra04a}, arcs and blocking sets~\cite{BKW05a} with the Kramer-Mesner-approach. In all these cases it is reasonable to investigate the limits of the refined Kramer-Mesner-method.

Our metaheuristic solver is a basic version of the tree search method Beam-search---there are several approaches to improve this method, for example by hybridizing it with another heuristic like ant colony optimization \cite{Blu05} or by combining it with a depth first search method \cite{FK05,ZH05}. Another promising approach could be the definition of a bounding function that weakens or replaces the proposed objective function by cutting feasible solutions which do not satisfy the bounds at some points in order to avoid computations of less promising solutions. 

\section*{Acknowledgement}

The authors would like to thank Patric \"Osterg\aa{}rd for his valuable comments to improve this manuscript.

\end{document}